\documentclass[]{amsart}

\newtheorem{theorem}{Theorem}[section]

\newtheorem{hypothesis}[theorem]{Unproven hypothesis}
\newtheorem{experiment}[theorem]{Experimental observation}

\theoremstyle{definition}

\theoremstyle{remark}

\numberwithin{equation}{section}

\usepackage[dvips]{graphics}

\title{Selfsimilarity and growth in Birkhoff sums for the golden rotation}
\date{May 31, 2010}
\author{Oliver Knill}
\email{knill@math.harvard.edu}
\address{Department of Mathematics, Harvard University, Cambridge, MA 02138,USA}
\author{Folkert Tangerman}
\email{tangerma@math.sunysb.edu}
\address{Department of Mathematics, Stony Brook University, Stony Brook, NY 11794,USA}

\begin{document}
\maketitle

\begin{abstract}
We study Birkhoff sums $S_k(\alpha) = \sum_{j=1}^k X_j(\alpha)$ with
$$ X_j(\alpha) = g(j\alpha) = \log| 2-2 \cos(2\pi j \alpha)| $$
at the golden mean rotation number $\alpha = (\sqrt{5}-1)/2$ with 
periodic approximants $p_n/q_n$. 
The summation of such quantities with logarithmic singularity is 
motivated by critical KAM phenomena. We relate the boundedness of $\log$ 
averaged Birkhoff sums $S_k/\log(k)$ and the convergence of $S_{q_n}(\alpha)$ with the 
existence of an experimentally established limit function
$f(x) = \lim_{n \to \infty} S_{[x q_n]}(p_{n+1}/q_{n+1})-S_{[x q_{n}]}(p_n/q_n)$ 
on $[0,1]$ which satisfies a functional equation $f(\alpha) + \alpha^2 f= \beta$ 
with a monotone function $\beta$. The limit $\lim_{n \to \infty} S_{q_n}(\alpha)$
can be expressed in terms of $f$. 
\end{abstract}

\section{Introduction}

Let $g$ be a periodic function of period $1$ and average $\int_0^1 g(x) dx=0$. 
For irrational $\alpha$, we consider the Birkhoff sums 
$S_k(g,\alpha)= \sum_{j=1}^k g(j \alpha)$.
Let $p_n/q_n$ be the sequence of continued fraction approximants to $\alpha$. 
It is known that if $g$ is real analytic and if $\alpha$ is sufficiently Diophantine, then $S_k$ stays bounded
and $S_{q_n}$ converges \cite{H}. Denjoy-Koksma theory \cite{CFS} which roots in work 
of Ostrowski, Hecke and Hardy-Littlewood in the twenties assures that if $g$ is of bounded variation, 
and $\alpha$ is of constant type, then $S_k(g,\alpha) \leq C \log(k)$ for some constant $C$ which only 
depends on $f$. If $g$ is continuous, then the boundedness of the sequence $S_k$ is by Gottschalk-Hedlund equivalent 
to the existence of a function $h$ such that $g(x) = h(x+\alpha)-h(x)$ \cite{TH,KH}.

Critical KAM phenomena \cite{KL,T} lead us to the study of
Birkhoff sums for the function 
\begin{equation}
   g(x)=\log(2-2 \cos(2\pi x)) =2 \log|2 \sin(\pi x)|=2\log|1-e^{2 \pi i x}|
\label{ourfunction}
\end{equation}
We assume throughout this paper that $\alpha$ is the golden mean $\alpha=(\sqrt{5}-1)/2$. 
The dynamical system $T(x) = x+\alpha \; {\rm mod} \; 1$ with this rotation number is the best 
understood and simplest aperiodic dynamical system. \\

While we only look at one orbit, one can sum along orbits of a more general deterministic random walk
$S_k(x,\alpha) = \sum_{j=1}^k X_j(x,\alpha)$, where
$X_j(x,\alpha) = g(x + j\alpha)$ are $L^1$ random variables over the probability space $(T^1,dx)$. 
In our case, the random variables $X_j$ have zero mean but nonzero quasi-periodic correlation. 
By Birkhoff's ergodic theorem for uniquely ergodic situations \cite{F}
we expect $S_k/k$ to converge to zero for every initial point not hitting the singularity $0$.
While we know $S_k/k \to 0$, we can look at the behavior of $\log$ averages $S_k/\log(k)$.
The choice of $\log(k)$ is due to the fact that for this particular Kronecker system and for functions 
$g$ of bounded variation we have $S_k = O(\log(k))$.  
It is natural to ask, to which extent limit theorems like the central limit theorem
apply, which are so well known when the random variables $X_k$ are independent. 
In the case of irrational rotations, 
the central limit theorem and the law of iterated logarithm typically use the same scaling function $\log(k)$
if $g$ is of bounded variation.  In our particular situation, where $g$ is given in (\ref{ourfunction}),
we see $\limsup_k S_k/\log(k) =2$ and $\liminf_k S_k/\log(k)=0$.
Experimentally, we measure that the distribution of $S_k/\log(k)$ to converge weakly to a 
stable distribution of compact support. 
One has to compare this with independent random variables of zero mean, where a scaling factor
$\sqrt{k}$ is needed to get a limiting distribution
and a rescaling by $\sqrt{2 k \log \log(k)}$ to have the sequence in a bounded interval almost surely. 
We will assume that our initial iteration point $x_0=0$ of the orbit agrees with the logarithmic singularity.
The positive orbit $x_n$ never hits the singularity $x_0$. \\

Birkhoff sums of functions with singularities over irrational rotation have been studied by
Hardy and Littlewood already \cite{HL}. They looked at $g(x) = \sin(x)^{-1}$ and showed that
the averaged partial Birkhoff sums $S_k/k$ stay uniformly bounded. 
In \cite{SU} an other non-integrable case $g(x) = (1-e^{ix})^{-1}$ is considered and a limiting
distribution for Birkhoff averages $S_k/k$ is established. Note that our case $g(x) = 2\log|1-e^{ix}|$ is 
integrable: $g \in L^1({\mathbb R})$ and we aim for a limiting distribution of $\log$ Birkhoff averages 
$S_k/\log(k)$ for a particular orbit. \\

Since $g$ has a logarithmic singularity at $x=0$, the function $g$ is not of bounded variation. 
When $\alpha$ is Diophantine of constant type, a modified Denjoy-Koksma argument shows that 
$|S_{q_n-1}(\alpha)|$ grows at most logarithmically so that $|S_k| \leq \log(k)^2$ 
(see \cite{KL}, where the idea there was to truncate the function and to use
standard Denjoy-Koksma for a truncated function). 
The orbit $j \alpha( \; {\rm mod} \; 1)$, $j=1, \dots,q_n-1$ remains away from the singularity $0$
and returns at the time $q_n$ just close enough so that $S_{q_n-1}$ grows logarithmically like  $\log(q_n)$
while $S_{q_n}$ converges.  \\

\begin{figure}
\scalebox{0.9}{\includegraphics{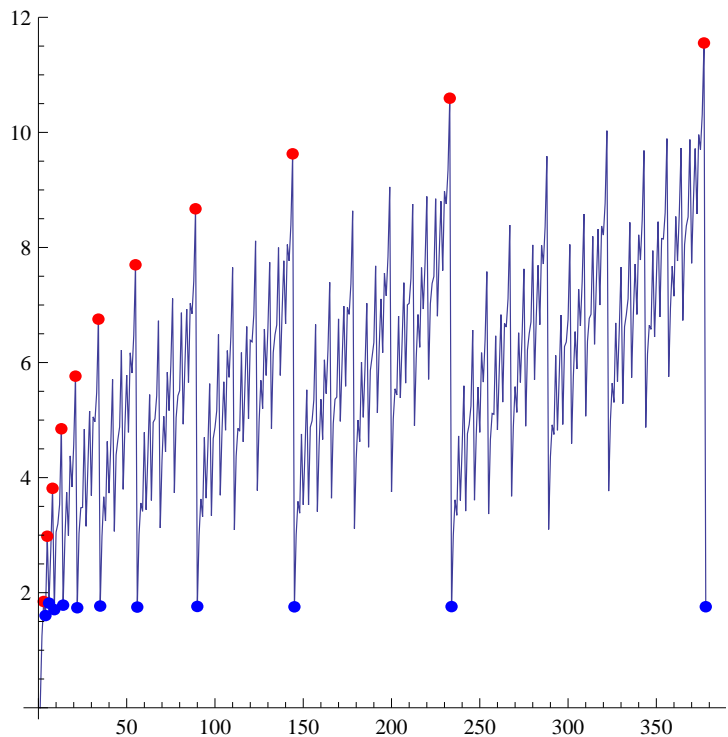}} \\
\caption{
The picture shows the sequence  $S_k(\alpha)$ between $k=1$ and $k=q_{14}$.
There are two envelopes: a lower one at around $1.75686 \dots$
corresponding to the converging subsequence $S_{q_k}(\alpha), k=1,\dots,14$, and an upper envelope 
that corresponds to the logarithmically diverging sequence $S_{q_n-1}(\alpha)$.
The values $S_{q_n}(\alpha)$ for $n=1,\dots ,14$ illustrate the existence of a positive limit.
We are able to compute this limit (see (\ref{numericalvalue})). 
}
\label{raw sequence}
\end{figure}

We make here the following three experimental observations: \\

\begin{center}\fbox{\parbox{13cm}{
\begin{experiment} 
$\;$  \\
\begin{enumerate} 
\item The sequence $S_{q_n}(\alpha)$ is \textbf{convergent} and has a finite \textbf{non-zero} limit. 
\item The sequence $S_k(\alpha)/\log(k)$ is {\bf bounded} and takes values in the interval $[0,2]$. 
\item The sequence $S_k(\alpha)/\log(k)$ has a {\bf symmetric limiting distribution} in $[0,2]$. 
\end{enumerate}
\end{experiment}
}}\end{center}

\vspace{2mm}
These go beyond \cite{KL}, where it was shown that $S_{q_n-1}(\alpha)/\log(q_n)$ is bounded leading to
$|S_k(\alpha)|\leq C (\log(k))^2$ for all $k$. 
We suggest a renormalization explanation for the first two observations
in Section~\ref{section3}.
Note that the convergence $1)$ implies the boundedness in $2)$ 
as indicated in the proof of Denjoy-Koksma. This step of the proof of 
Denjoy-Koksma's theorem does not need that $g$ is of bounded variation.  \\

Birkhoff sums over an irrational rotation are estimated as follows if the function has bounded variation.
If $p/q$ is the periodic continued fraction approximation of 
$\alpha$, then $|S_q| \leq C$ for all $q$ implies $|S_k| \leq  C_1 k^{1-1/r} \log(k)$ with a $k$ independent
constant $C$, if $\alpha$ is of Diophantine type $r$ (see i.e.\cite{J,CFS}). 
In the golden mean case, where $r=1$, this gives $|S_k| \leq  K \log(k)$ for some other constant $K$. 

\begin{figure}
\scalebox{0.7}{\includegraphics{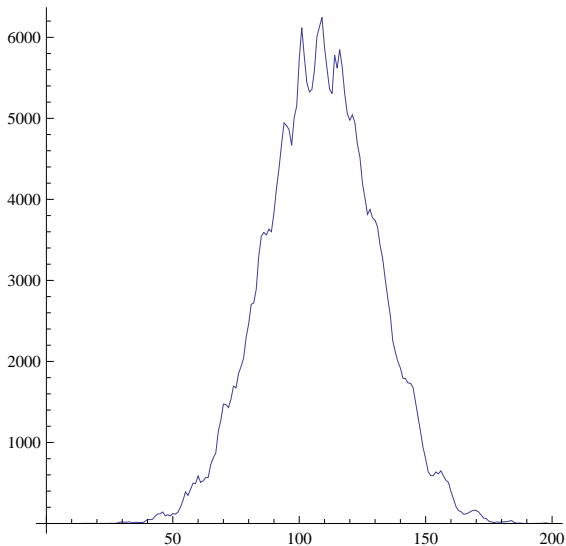}} \\
\caption{
The distribution $\mu$ of the sequence $S_k(\alpha)/\log(k)$,
$k=1,\dots,q_{28}$ is sampled at 200 uniformly spaced
bins covering the interval $[0,2]$. Mathematically, the distribution needs
to be understood as the functional
$\phi \to \lim_{n \to \infty} \frac{1}{q_n} \sum_{k=1}^{q_n} \phi(S_k(\alpha)/\log(k)) = \int \phi \; d\mu$.
The closest distribution we can think of is a variant of a Bernoulli convolution \cite{BBG}. 
At least, it has similar features. 
}
\label{distribution}
\end{figure}

\section{Origin of the Birkhoff sums}

Birkhoff sums have been studied in the context of number theory, probability 
theory and dynamical system theory. Birkhoff's ergodic theorem $S_k/k \to \int_0^1 g(x) \; dx$ 
relates Birkhoff averages $S_k/k$ with the average of $g$. If $g$ is continuous, this convergence is
independent of the initial point if the dynamical system is an irrational rotation \cite{F}. 
If the mean $\int g(x) \; dx=0$, one can ask how fast $S_k$ grows. In the case of 
rigid rotation, the answer depends on the regularity of $g$ and arithmetic properties of the rotation 
number. In the case of an irrational rotation we have $S_k = O(1)$ for Diophantine 
$\alpha$ and real analytic $g$. For functions $g$ of bounded variation, Denjoy-Koksma theory allows 
estimates $S_k  = O(k^r)$ for $0 \leq r<1$ and of the form $S_k = O(\log(k))$ for $\alpha$ of 
constant type. There are also lower bounds on the growth rate: Herman \cite{TH}, using it as a tool for
higher dimensional considerations,  gave example of Birkhoff sums, 
where $\limsup S_k \geq  \sqrt{k}$. Bi\`evre and Forni have for any $\epsilon>0$ examples which achieve 
$\limsup S_k \geq k^{1-\epsilon}$ \cite{BF}. We look here at an $L^1$-integrable case of unbounded variation. \\

\begin{center} \begin{tabular}{lll}  \hline
Hecke                      &  $g(x) = \lfloor x \rfloor = x-[x]$  & \cite{hecke1922} 1921 \\
Hardy-Littlewood           &  $g(x) = \sin(x)^{-1}$               & \cite{HL} 1928    \\
Bryant-Reznick-Serbinowska  & $g(x) = (-1)^{\lfloor x \rfloor}$   & \cite{BS} 2006    \\
Sinai-Ulcigrai             &  $g(x) = 1/(1-\exp(i x))$            & \cite{SU} 2008    \\
This paper                 &  $g(x) = \log|1-\exp(i x)|$          & 2010              \\ \hline 
\end{tabular} \end{center} 

\vspace{5mm}
The authors of the present paper were led to the study of these sums from different vantage points.
The first author (K) was motivated by the study of determinants of truncated versions of the
diagonalization of the linear operator,
\begin{equation} h \to Lh =  h(x+\alpha)+h(x-\alpha)-2h(x) \label{operator} \end{equation}
on $L^2(\mathbb{T}^1)$, where $\mathbb{T}^1$ is the unit circle (see \cite{KL}). Perturbations of 
the operator $L$ appear in the study of invariant circles in twist maps and 
the growth rate of determinants of truncated Fourier transformed matrices are used to estimate 
Green functions \cite{B}. \\
The second author (T) was motivated by the solution of the linear conjugacy relation:
\begin{equation} m(\lambda z)+m(\lambda^{-1}z)-2m(z)=zm(z) \label{conjugacy equation} \end{equation}
where $\lambda = e^{2\pi i \alpha}$ and $m=\sum_{k=1}^{\infty} m_k z^k$ is a holomorphic function with 
$m(0)=0$.
The coefficients $m_k$ are positive and $m_k=e^{S_k(\alpha)}$.
The logarithmic growth of the coefficients, $S_k(\alpha)$
shows that the power series for $m$ is convergent in the interior of the unit disk.
The distribution properties of the sequence $S_k$ shows that this sequence does not converge at $z=1$.
It it is not known if $m$ has a meromorphic extension to the entire complex plane.
Equation~(\ref{conjugacy equation}) occurs in various KAM contexts. This equation occurs
also a first order solution occurring in the study of the boundary of KAM rings in
conservative systems, see \cite{SV,T}. \\

An other motivation connects the present study with the field of holomorphic dynamics: 
consider the nonlinear complex dynamical system in $C^2$:
$$     T(z,w) = (c z , w (1-z) ) \; ,  $$
where $c = \exp(2\pi i \alpha)$ and $\alpha$ is the golden mean. 
This is one of the simplest quadratic systems which can be written down in $\mathbb{C}^2$. How
does the orbit behave on the invariant cylinder $\{ |z|=1 \} \times \mathbb{C}$ starting at $(c,1)$? We have
$$  T^n(z,w) = (z_n,w_n) = ( c^n z,  w (1-z) (1-c z) (1-c^2 z) \dots (1-c^{n-1} z) )  $$
and
$$ \log|w_n| =  \sum_{k=1}^n \log| 1-e^{2\pi i k \alpha} | = \frac{S_{n}}{2} $$
for $(w_0,z_0)=(c,1)$ because $2 \log|1-e^{ix}| =  \log(2-2 \cos(x))$.
The study of the global behavior of the holomorphic map $T$
in $\mathbb{C}^2$ boils down to the Birkhoff sum over the golden circle on a subset because
for $r=|z|<1$, where $g_r(x) = \log|1-r e^{ix}|$ is real analytic and the Birkhoff sum converges by Gottschalk-Hedlund.
It follows that for $r<1$ the orbits have the graph of a function $A: \{ |z|=r \; \} \to \mathbb{C}$ as 
an attractor. For $r=|z|>1$, we have $|w_n| \to \infty$. So, all the nontrivial dynamics of the quadratic 
map happens on the subset $\{ |z|=1 \; \} \times \mathbb{C}$. \\

There is also a relation with the Dirichlet series $\sum_{k=1}^n a_k/\log(k)^s$ with 
$a_k = \log(2-2 \cos(k \alpha))$ which has the 
abscissa of convergence $\limsup_n \log|S_k|/\log\log(k)$. The experiments indicate that the 
abscissa of convergence is $1$.  

\section{A renormalization picture}
\label{section3}
In order to shed light on the main observation, we rescale the summation interval $[1,q_n]$ 
to the unit interval $[0,1]$. More precisely, let $p_n/q_n$ be the periodic approximants 
of $\alpha$ and let $[t]$ be the largest integer smaller than $t$. We consider the 
piecewise constant function $S_{[xq_n]}(\alpha)$ which has discontinuities at points
$x=\dfrac{p_n}{q_n}$ for which $[xq_n]$ is an integer. Rescaling the sequence 
this way, we can talk about functions on the unit interval $[0,1]$ and can treat all 
the sequences in the same space of functions. \\

By looking at the convergence of Fourier coefficients, 
we observed a weak limit $s(x)$ of $s_n(x)=S_{[x q_n]}(\alpha)/\log([x q_n])$.

\begin{figure}
\scalebox{1.2}{\includegraphics{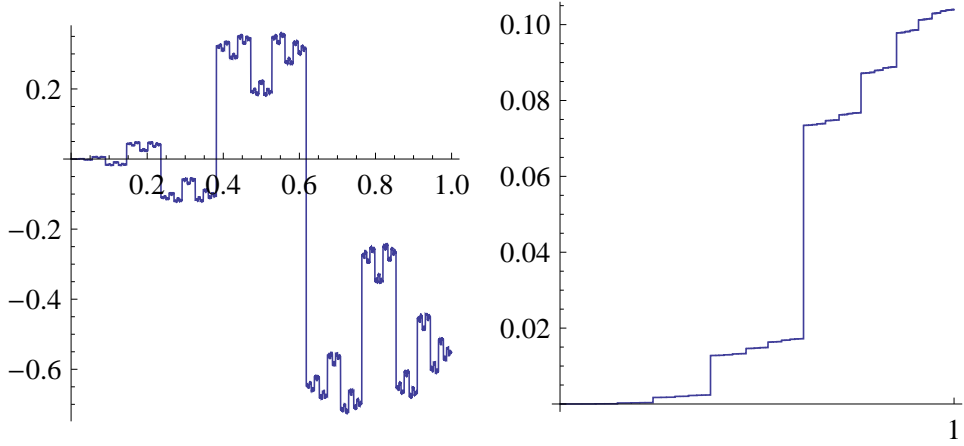}}
\caption{
\label{renormalization}
The "skyline" function $f_{24}(x)$ and the "stair" function $\beta_{24}(x)$. 
The function $\beta$ appears monotone with $\beta(1-)=0.104 \dots $.
}
\end{figure}

In order to study the limit better and see pointwise convergence we look at the sequences
$$  S_k( \frac{p_{n+1}}{q_{n+1}}) -  S_k(\frac{p_{n}}{q_{n}}) $$
with $k \in \{ 1, \dots , q_n \; \}$. In the rescaled picture, we have functions
$$  f_n(x) = S_{[x q_n]}( \frac{p_{n+1}}{q_{n+1}}) -  S_{[x q_n]}(\frac{p_{n}}{q_{n}}) $$
on the interval $[0,1]$. These functions $f_n$ appear to converge (see Figure~(\ref{renormalization})). \\

While some self similarity is present for the original functions $s_n(x)$, where they only manifest weakly,
the new functions $f_n$ have more regularity, appear to converge pointwise almost everywhere and appear to have a
self-similarity property: \\

\fbox{\parbox{12cm}{
\begin{hypothesis}
$\;$ \\
The functions $f_n$ converge pointwise almost everywhere to a function $f$ satisfying  $|f(x)| \leq x$ and 
\begin{equation} f(\alpha x) + \alpha^2 f(x) =  \beta(x)  \; , \label{selfsimilarity} \end{equation}
where $\beta(x)$ a monotone increasing function which is zero for $x=0$.
\label{hypothesis}
\end{hypothesis}
}}

\begin{figure}
\scalebox{1.2}{\includegraphics{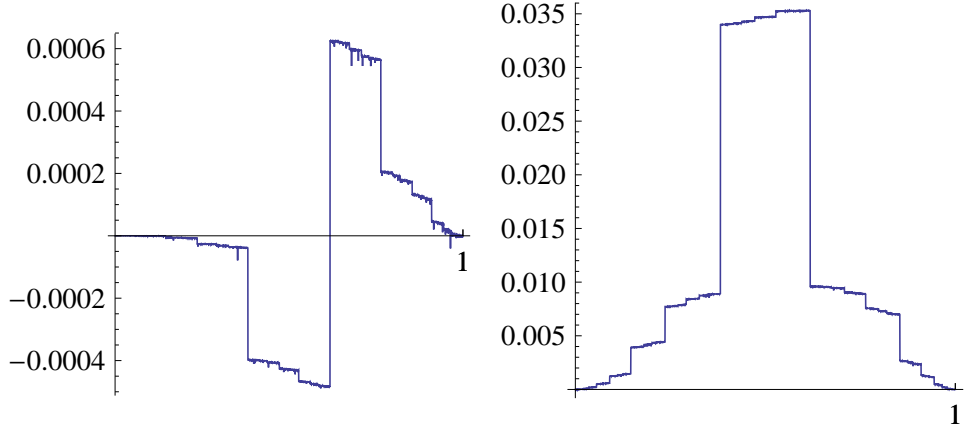}}
\caption{
\label{betafine}
These graphs show a comparison of the entire graph with the lower and upper part of $\beta$.
To the left we see the function $\beta-\beta_1$, to the right, we see the function $\beta-\beta_2$.
}
\end{figure}

\vspace{4mm}

The function $\beta$ appears to be the anti derivative of a positive measure $\mu$ so that $\beta(x) = \int_0^x \; d\mu(x)$. 
Discontinuities occur along the forward orbit 
$x = n\alpha$, $n=1,2, \dots$. We observe that the graph $\beta(x)$ looks like a "devil stair case" and is close 
to self similar: the graph of $\beta$ on $[0,\alpha]$ matches the graph of $\beta$ on $[0,1]$ 
closely: define $a=\beta(\alpha-), b=\beta(\alpha+),c=\beta(1-)$ and 
the functions $\beta_1(x) = \beta(\alpha x) (c/a)$ and
$\beta_2(x) = (\beta(\alpha+ (1-\alpha)x)-b) c/(c-b)$. We observe that they are both close to $\beta$.
The comparison of the affine scaled graphs of $\beta$ are seen in Figure~(\ref{betafine}). \\

The identity:
$$ S_k(\alpha)-S_k(\dfrac{p_N}{q_N})= \sum_{n=N}^{\infty}S_k(\dfrac{p_{n+1}}{q_{n+1}})-S_k(\dfrac{p_n}{q_n}) $$
suggests the introduction of a sequence of functions:
\begin{equation} h_n(x) = S_{[x q_n]}(\alpha) - S_{[x q_n]}(\frac{p_n}{q_n}) \label{hn} \end{equation}
\begin{equation} f_{n,m}(x) = S_{[x q_n]}( \frac{p_{n+m}}{q_{n+m}}) - S_{[x q_n]}( \frac{p_n}{q_n}) \label{fnm} \end{equation}
By definition, $f_n = f_{n,1}$ and $h_n = f_{n,\infty}$, where we understand $p_{\infty}/q_{\infty} = \alpha$. \\

\section{Consequences of the hypothesis}

The following statements do follow from the yet unproven hypothesis. \\

\vspace{0.4cm}

{\bf 1. The functions $h_n$ converge to a function $h(x)$ determined by $f$.} \\ 

Proof.  
For every $j>0$ and a set $Y_n$ of $x$ of measure $|Y_n|$ converging to $1$, we have
$f_{n,j+1}(x) - f_{n,j}(x) = f_{n}(\alpha^j x)$.  Using $f_n=f_{n,1}$, one obtains that
$$ f_{n,m}(x) - f_n(x) = \sum_{j=0}^{m-1} f_{n}(\alpha^j x)  \;  $$
and $h_n=f_{n,\infty}(x) - f_n(x)  = \sum_{j=0}^{\infty} f_{n-1}(\alpha^j x)$. 
Since by the hypothesis, $f_n$ converges pointwise almost everywhere to a function $f$,
we know that $|f_n(x)| \leq |x|$ for large enough $n$. For fixed $n$, 
the sup-norm of $f_n(\alpha^j)$ decays at an exponential rate 
$\alpha^j$ and $h_n$ is bounded and has a limit 
\begin{equation}
\label{fandh}
 h(x)=\sum_{j=0}^{\infty} f(\alpha^jx) \; . 
\end{equation}
The graph of the function $h$ looks similar than the graph of $f$. \\

While $f$ was the difference between Birkhoff sums with periodic rotation numbers, the function
$h$ shows the difference between Birkhoff sums of irrational and rational rotation numbers in the limit. 
The function 
$$  \gamma(x) = h(\alpha x) + \alpha^2 h(x) =\sum_{j=1}^{\infty} \beta(\alpha^j x) $$ 
is a monotone function too and looks like $\beta$. \\

{\bf 2. Derivation of the experimental observation: the limit along subsequences are known. } \\

By definition,  $h_n(1^-) = S_{q_n-1}(\alpha) - S_{q_n-1}(p_n/q_n)$.
Since 
$$   S_{q_n}(\alpha) - S_{q_n-1}(\alpha) + 2\log(q_n) \to \log(\frac{4\pi^2}{5})=2.06632\dots \;  $$
for $n \to \infty$ (see (\ref{fourpisquarefifth}) below),
and 
$$  S_{q_n-1}(\frac{p_n}{q_n}) = 2 \log(q_n) $$ 
noted in \cite{KL}, we have
\begin{equation}
\label{numericalvalue}
S_{q_n}(\alpha) \to \log(\frac{4\pi^2}{5}) +h(1^-) = 1.75687 \dots \; 
\end{equation}
where $h(1^-) = \lim_{x \nearrow 1} h(x)= -0.30945 \dots$. We see that the limiting value
of $S_{q_n}(\alpha)$ can be read off from the function $h$ which is related to $f$
by (\ref{fandh}).  \\

{\bf 3. Derivation of the experimental observation: The sequence $S_k(\alpha)/\log(k)$ is bounded.} \\

Proof. This part follows from the second part of the proof of Denjoy-Koksma in the cases 
when $\alpha$ is Diophantine of bounded type. 
From $S_{q_n}(\alpha) \leq M$ follows $S_k(\alpha) \leq M_2 \log(k)$ for all $k>0$ 
with an other constant $M_2$ independent of $k$.  \\

{\bf 4. Derivation of the experimental observation: $S_k(\alpha)/\log(k)$ have accumulation points in $[0,2]$. } \\

The statement
$$ \lim_{n \to \infty} S_{q_{n}}(\alpha)/\log(q_n)\,=0 $$
follows immediately from the convergence of $S_{q_n}(\alpha)$.  The statement
$$ \lim_{n \to \infty} S_{q_{n}-1}(\alpha)/\log(q_n)\,=2 $$
follows from 
$$ S_{q_{n}}(\alpha)  - S_{q_{n}-1}(\alpha) = X_{q_n}(\alpha) = \log|2-2 \cos(2\pi q_n \alpha) | $$
and 
$$  \lim_{n \to \infty} q_n^2 ( 2-2\cos(2\pi q_n \alpha) )  = 
    \lim_{n \to \infty} q_n^2 ( 2-2\cos(2\pi q_n \alpha-p_n) ) = \frac{4\pi^2}{5} \;  $$
which implies 
\begin{equation}
\label{fourpisquarefifth}
 X_{q_n}(\alpha) + 2 \log(q_n) \to \log(\frac{4\pi^2}{5}) \; . 
\end{equation}
We especially see that the normalized sequence $\{S_k/\log(q_n) \}_{k=1}^{q_n}$ jumps by $-2 + O(1)/\log(q_n)$
when $k=q_n$. 

\section{Conclusions and Open Questions} 

The particular Birkhoff sum for the golden mean rotation shows interesting patterns.
It is a situation where Denjoy-Koksma falls short of explaining the experimentally measured growth rate. 
The nature of the limiting distribution $\rho$ of $S_k(\alpha)/\log(k)$ is not yet settled. 
Taking test functions $\phi(x)_m=e^{i m x}$ gives the Fourier coefficients 
$$  \hat{\rho}_m = \lim_{n \to \infty} \frac{1}{q_n} \sum_{k=1}^{q_n} \phi_m(\frac{S_k(\alpha)}{\log(q_n)}) $$
of this measure. Understanding the still unproven hypothesis could answer the question whether the limiting
distribution exists and if yes, how these coefficients decay. In this paper, 
we have traced the convergence and boundedness questions to the existence of limiting functions
given in hypothesis~(\ref{hypothesis}). \\

\begin{figure}
\scalebox{1.2}{\includegraphics{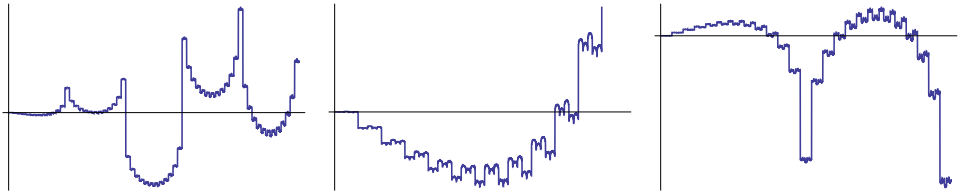}} \\
\caption{
\label{otherrotationnumber}
For $\alpha=\sqrt{41}-6$ with continued fraction expansion $[2,2,12,2,2,12,\dots]$ of period $L=3$,
each of the three function sequences 
$\{ f_{3n} \; \}_{n \in {\mathbb{N}}},
 \{ f_{3n+1} \; \}_{n \in {\mathbb{N}}}, 
 \{ f_{3n+2} \; \}_{n \in {\mathbb{N}}}$ 
seem to converge. The picture shows these three graphs in the case $n=3$. If $\alpha$ has a continued
fraction expansion of period $L$, we see that $f_{Ln}$ and $S_{q_{Ln}}$ do converge. Similar statements as in
hypothesis~(\ref{hypothesis}) seem to hold: for $L=2$ with examples like $\alpha=(\sqrt{3}-1)/2=[2,1,2,1,\dots]$,
we see that $f_{2n}$ converges to a function $f$ and $f_{2n+1}$ converges to a function $\tilde{f}$ and that
$\beta(x) = f(\alpha x) + 2 \alpha^2 \tilde{f}(x)$ is monotone. For
$\alpha=[p,q,p,q,\dots]$ we observe a monotone function 
$\beta(x) = q f(\alpha x) + p\alpha^2 \tilde{f}(x)$. 
}
\label{raw sequence}
\end{figure}

1. {\bf Other rotation numbers}. 
For $(\sqrt{2}-1)/2=[4,1,4,1,4,1,\dots]$ for example, where the continued fraction expansion has
period $2$, there appears to be a function $f$ such that $f_{2n}$ converges to $f$ while $f_{2n+1}$ converge to $-f$.
For $\sqrt{41}-6=[2,2,12,2,2,12,\dots]$, where the continued fraction expansion has period $3$, the sequences 
$f_{3n},f_{3n+1},f_{3n+2}$ seem to converge (see Figure~\ref{otherrotationnumber}). 
This suggests that for quadratic irrationals, there is presumably a similar story
depending on the structure of the continued fraction expansion. The golden mean $(\sqrt{5}-1)/2=[1,1,1,1,\dots]$ is the simplest.
The silver ratio $\sqrt{2}-1=[2,2,2, \dots]$ similarly as the golden ratio shows a limit for $f_n$. 
For other Diophantine rotation numbers that satisfy a Diophantine or Brjuno condition, 
the growth is expected to match the estimates given by Denjoy-Koksma in the case of bounded variation. \\

2. {\bf Different starting point}. Starting at a different point $x_0$ changes the story.
This is not surprising since it is no more the $q_n$ which lead to the close 
encounters with the logarithmic singularity but by a theorem of Chebychev~\cite{K}, 
there are integers $k<q_n, l$ with  % THeorem 24}
$|x_0+k \alpha -l|<3/q_n$. For the self-similarity structure we have seen here, it looks as if it is
important to have the initial point $x_0$ at the logarithmic singularity.  \\

3. {\bf Other functions}. If $g(x) = \log(2-2\cos(x))$ is replaced with an other function 
$\log(r(x))$ with one logarithmic singularity at $0$, the story is just distorted:  
Assume $r(x)$ is a trigonometric polynomial with a single root at $0$ having the property that $r''(0)=2$.
Then $g(x) = \log(r(x))$ shows similar growth rates. The reason is that 
$g(r(x)) - \log(2-2 \cos(x)) = \log(r(x)/(2-2 \cos(x)))$ is now smooth
and has bounded variation. Adding a function of bounded variation to $g$ does not change the behavior
because of Denjoy-Koksma theory. It is important however that we start the orbit at the critical point. \\

\begin{figure}
\scalebox{0.8}{\includegraphics{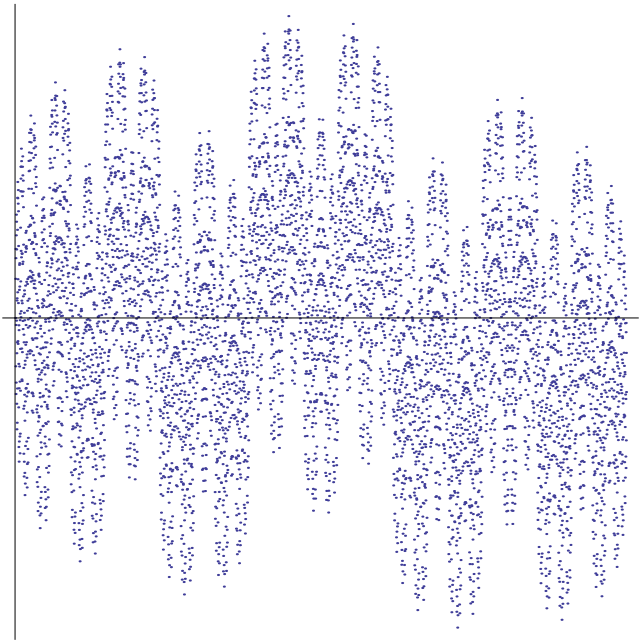}}
\label{hecke}
\caption{
We see the Birkhoff sum $\{ S_k(\alpha), \; 1 \leq k \leq q_{20} \; \}$ in Hecke's case \cite{hecke1922}, 
where $g(x) = x-[x]$. Here, $S_{q_{2n}}(\alpha),S_{q_{2n+1}}(\alpha)$ converge and the 
functions $f,h$ are explicitly known: The functions $h_n(x) =  S_{[x q_n]}(\alpha) - S_{[x q_n]}(p_n/q_n)$  have
the property that $h_{2n}(x)$ converges pointwise to the function $-c x^2$ with $c = \pi/\sqrt{5}$ and
$h_{2n+1}(x)$ converges pointwise to $c x^2$. From $h=\sum_{n=0}^{\infty} f(\alpha^k x)$
and the asymptotic relation $f_n(\alpha x) = -\alpha^2 f_{n-1}(x)$ with
$f_n(x) = S_{[x q_n]}(p_{n+1}/q_{n+1})-S_{[x q_n]}(p_n/q_n)$ we get 
$h(x)= \sum_{n=0}^{\infty} (-\alpha^2)^n f(x) = f(x)/(1+\alpha^2)$ so that $f=\lim_{n \to \infty} f_{2n}$
is a multiple of $h=\lim_{n \to \infty} h_{2n}$. }
\end{figure}

4. {\bf Selfsimilarity in Hecke's example}. Birkhoff sums for bounded $g$ have been studied first by number theorists
like Hardy,Littlewood or Hecke. In the case of a golden rotation, there is a similar structure also but the story is different
in that the limiting functions $f$ and $h$ are smooth. 
Figure~(\ref{hecke}) shows the Birkhoff sum in the situation studied by Hecke \cite{hecke1922} where $g(x)=x-[x]$ is
piecewise smooth. It is historically the first example studied for irrational rotation numbers.  
The functions $f = \lim_n f_{2n}$ and $h=\lim_n h_{2n}$ are explicitly known quadratic function
in the case when the rotation number $\alpha$ is the golden mean.  \\

5. {\bf Replacing the irrational rotation by a smooth circle map}. 
A similar deformation happens if one replaces the irrational golden rotation with a 
smooth circle map $T$ with the same rotation as long the conjugated orbit starts at the critical point of $g$.
The reason is that any smooth interval map $T$ with 
Diophantine rotation number is conjugated to to an irrational rotation: $T^n(x) = S^{-1}(S(x)+n \alpha))$ 
and so $g(T^nx) = g(S^{-1}(S(x) + n \alpha))$. Therefore, if the starting point is $x=S^{-1}(0)$, then 
$g(T^nx) = h(n \alpha)$ with $h(x) = g(S^{-1}(x))$ so that $h$ and $g$ have the same logarithmic singularity
and $g(x) =  h(S(x))$. We see that changing to a smooth circle map has the same effect as replacing $g$ with 
an other function with the same logarithmic singularity and changing the initial point.  If the initial point
is the same, we see the same behavior but distorted by a Denjoy-Koksma correction for a continuous function $g$.  \\

\begin{figure}
\scalebox{1.0}{\includegraphics{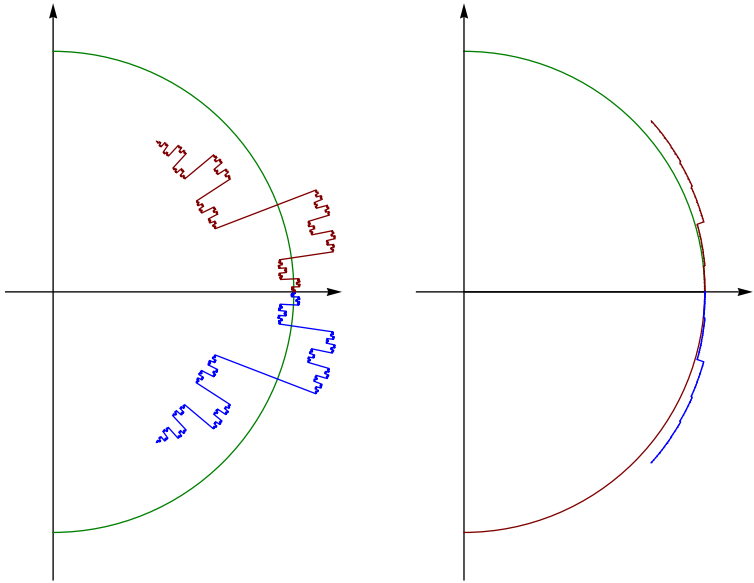}}
\caption{
\label{complex}
The figure to the left shows simultaneously the images of the rational functions 
$F_{15}$ and $F_{16}$ from $[0,1]$ into the complex plane.
It appears that the images of $F_{2n}$ and $F_{2n+1}$ converge pointwise
to complex conjugates. We have $\log|F_n|=f_n$. To the right we see the image of
the function $B_{15}$ and $B_{16}$ in the complex plane.
It appears as if the images of $B_{2n}$ and $B_{2n+1}$ converge pointwise and that the limits
are complex conjugate. We have $\log|B_n|=\beta_n$. }
\end{figure}

6. {\bf Replacing the irrational rotation by a toral map}. 
Consider $T(x,y)= (x+\alpha,y+\beta)$ and $g(x,y) = g(2\pi x y)$ and $\gamma=\alpha \beta$, we get Birkhoff sums
$$  S_k = \sum_{n=1}^k g(n^2 \gamma) \; .  $$
Illustrations of these Birkhoff sums are called {\bf curlicues} see \cite{BG,S2} where $g(x) =  \exp(2\pi i x)$
leads to Weyl sums $S_m=\sum_{k=1}^m \exp(2\pi i p(n))$ where $p$ is a polynomial.  A normalization $S_m/\sqrt{m}$
appears to lead to a limiting distribution \cite{S2}. \\

\begin{figure}
\scalebox{0.8}{\includegraphics{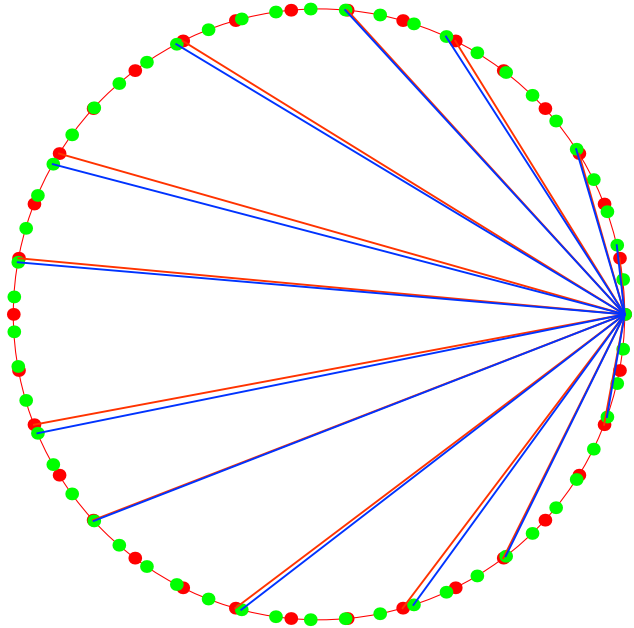}}
\caption{
\label{polygon}
Geometric interpretation of the absolute value of the complex function $F_n$ as a product of length quotients
of regular polygons with Fibonacci order. The picture illustrates the value $|F_9(0.4)|$. }
\end{figure}

7. {\bf Complex case}. Because of Equation~(\ref{ourfunction})
the Birkhoff sums originate in products
$$ P_n(\alpha)  = \prod_{k=1}^n (1-e^{2\pi i k \alpha})   $$
which satisfy 
$$ \log| P_n(\alpha)|  = \frac{S_n}{2} \; . $$
The functions
$$  F_n(x) = \prod_{k=1}^{[x q_n]} \frac{ 1-e^{2\pi i k \frac{p_{n+1}}{q_{n+1}}}}{ 1-e^{2\pi i k \frac{p_n}{q_n}}},   \; $$
on $[0,1]$ satisfy $f_n=\log|F_n|$. Pointwise convergence of $F_n$ appears to occur modulo complex conjugation.
For odd $n$, the image of $F_n$ is in the upper half plane, for even $n$ 
it is in the lower half plane. The functions
$$  B_n(x) = F_n(\alpha x) \cdot F_n(x)^{\alpha^2} $$
which satisfies $\log|B(x)| = \beta(x)$ appear to converge modulo complex conjugation. 
Establishing the convergence hypothesis in the complex would of course establish the hypothesis 
in the real case. See Figure~(\ref{complex}). \\

To establish the first experimental observation (and so the second), we would 
only have to show the existence of the limit
$$  \lim_{n \to \infty} \prod_{k=1}^{q_n-1} \frac{ 1-e^{2\pi i k \alpha}}{ 1-e^{2\pi i k \frac{p_n}{q_n}}} \; ,   \; $$
where $p_n/q_n$ are the Fibonacci quotients converging to the golden mean $\alpha$.  \\

8. {\bf Geometric interpretation}. Our observations have an elementary geometric interpretation:
take the Fibonacci sequence $\{q_n\}_{n=1}^{\infty} = (1,1,2,3,5,8,13,21, \dots)$ 
and consider the regular $q_{n+1}$-gon and $q_n$-gons in the unit circle. The number
$$ |P_{q_{n}-1}(\frac{q_{n-1}}{q_n})| $$ 
is the product of the lengths of all the diagonals in the
$q_n$-gon which start from one point. This product is $q_n$ \cite{KL}. \\

The value of $|F_n(x)|$ is a product of ratios between lengths of 
a subset of all diagonals in the $q_{n+1}$-gon and the $q_n$-gon. See Figure~(\ref{polygon}).

\vfill
\pagebreak

\bibliographystyle{plain}

\vspace{1cm}

\end{document}